\documentclass[11pt,reqno]{amsart}
\usepackage{amsmath,amssymb,amsthm,mathrsfs,dsfont}
\usepackage{CJK}

\newtheorem{thm}{Theorem}[section]
\newtheorem{lem}{Lemma}[section]
\newtheorem{defn}{Definition}[section]
\newtheorem{prop}{Proposition}[section]
\newtheorem{coro}{Corollary}[section]\numberwithin{equation}{section}
\newtheorem{rmk}{Remark}[section]

\def\pf{{\textit {Proof:} }}

\newcommand{\mysection}[1]{\section{#1}\setcounter{equation}{0}}

\newfont{\bb}{msbm10 at 12pt}

\newcommand{\bal}{\begin{aligned}}      \newcommand{\eal}{\end{aligned}}
\newcommand{\ba}{\begin{array}}      \newcommand{\ea}{\end{array}}
\newcommand{\bc}{\begin{center}}     \newcommand{\ec}{\end{center}}
\newcommand{\be}{\begin{enumerate}}  \newcommand{\ee}{\end{enumerate}}
\newcommand{\beq}{\begin{eqnarray}}  \newcommand{\eeq}{\end{eqnarray}}
\newcommand{\beQ}{\begin{eqnarray*}} \newcommand{\eeQ}{\end{eqnarray*}}
\newcommand{\bi}{\begin{itemize}}    \newcommand{\ei}{\end{itemize}}
\newcommand{\bt}{\begin{tabular}}    \newcommand{\et}{\end{tabular}}
\newcommand{\bdm}{\begin{displaymath}} \newcommand{\edm}{\end{displaymath}}

\def\qed{\hfill{Q.E.D.}\smallskip}

\newcommand{\ls}{\setlength{\baselineskip}{12pt}
                 \setlength{\parskip}{3mm}}

\title{2-Dimensional Combinatorial Calabi Flow in Hyperbolic Background Geometry}
\author{Huabin Ge, Xu Xu}
\address[Huabin Ge]{School of Mathematical Sciences, Peking University, Beijing 100871, PR China}
\email{gehuabin@pku.edu.cn}

\address[Xu Xu]{School of Mathematics and Statistics, Wuhan University, Wuhan 430072, PR China}
\email{xuxu2@whu.edu.cn}
\date{}

\begin{document}
\maketitle

\begin{abstract}
For triangulated surfaces locally embedded in the standard hyperbolic space, we introduce combinatorial Calabi flow as the negative gradient flow of combinatorial Calabi energy. We prove that the flow produces solutions which converge to ZCCP-metric (zero curvature circle packing metric) if the initial energy is small enough.  Assuming the curvature has a uniform upper bound less than $2\pi$, we prove that combinatorial Calabi flow exists for all time.
Moreover, it converges to ZCCP-metric if and only if ZCCP-metric exists.
\end{abstract}

\mysection{Introduction}\ls

Consider a compact surface $X$ with a triangulation $T=\{V,E,F\}$, where the symbols $V,E,F$ represent the set of vertices, edges and faces respectively. A positive function $r:V\rightarrow (0,+\infty)$ defined on the vertices  is called a circle packing metric and
a function $\Phi: E\rightarrow [0, \pi/2]$  is called a weight on the triangulation.
Throughout this paper, a function defined on vertices is regarded as a column vector and $N=V^{\#}$ is used to denote the number of vertices. Moreover, all vertices, marked by $v_{1},...,v_{N}$, are supposed to be ordered one by one and we often write $i$ instead of $v_i$. Thus we may think of circle packing metrics as points in $\mathds{R}^N_{>0}$, $N$ times of Cartesian product of $(0,\infty)$.

Given $(X,T,\Phi)$, every circle packing metric $r$ determines a piecewise linear metric on $X$ by attaching edge $e_{ij}$ a length $l_{ij}=\sqrt{r^2_i+r^2_j+2r_ir_jcos(\Phi_{ij})}.$ This length structure makes each triangle in $T$ isometric to an Euclidean triangle. Furthermore, the triangulated surface $(X,T)$ is composed by gluing many Euclidean triangles coherently. This case is called Euclidean background geometry in \cite{CL1}, where the cases of hyperbolic background geometry and spherical background geometry are also studied.
In these cases, the length $l_{ij}$ of the edge is determined by the hyperbolic and spherical cosine law respectively,
$$\cosh l_{ij}=\cosh r_i\cosh r_j+\sinh r_i\sinh r_j\cos \Phi_{ij},\;\; \text{for}\ \ \mathds{H}^2,$$
 $$\cos l_{ij}=\cos r_i\cos r_j-\sin r_i\sin r_j\cos \Phi_{ij},\;\;\text{for}\ \ \mathbb{S}^2.$$

Inspired by \cite{CL1} and \cite{Ge1}, we study the negative gradient flow of discrete Calabi energy on $X$ in hyperbolic background geometry.
The flow is called combinatorial Calabi flow \cite{Ge1}.
We use $(X,T,\Phi,\mathds{H}^2)$ to denote the space we want to study in the following
in hyperbolic background geometry, where $X$ is a closed surface, $T$ is a fixed triangulation of $X$,
$\Phi$ is a fixed weight function defined on edges, and $\mathds{H}^2$ represents hyperbolic background geometry .
It was proved by Thurston \cite{T1} that, whenever $\{i,j,k\}\in F$, these three positive numbers $l_{ij}, l_{ik}, l_{jk}$ satisfy the triangle inequalities. Thus the combinatorial triangle $\{i,j,k\}$ with lengths $l_{ij}, l_{ik}, l_{jk}$ form a hyperbolic triangle in $\mathds{H}^2$.

For the hyperbolic triangle $\triangle v_iv_jv_k$, the inner angle of this triangle at $v_i$ is denoted by $\theta_i^{jk}$,
and the combinatorial curvature $K_i$ at $v_i$ is defined as
$$K_i=2\pi-\sum_{\{i,j,k\}\in F}\theta_i^{jk}.$$
Notice that $\theta_i^{jk}$ can be calculated by hyperbolic cosine law,
thus $\theta_i^{jk}$ and $K_i$ are elementary functions of the circle packing metric $r$.

For the hyperbolic triangle $\triangle v_iv_jv_k$, we have
$\theta_i^{jk}+\theta_j^{ik}+\theta_k^{ij}=\pi-Area(\triangle v_iv_jv_k)$.
Using this formula, we have the combinatorial Gauss-Bonnet formula \cite {CL1}
\begin{equation}\label{Gauss-Bonnet formula}
\sum_{i=1}^NK_i=2\pi \chi(X)+Area(X).
\end{equation}

\section{Combinatorial Calabi Flow on Surfaces}\ls

\subsection{The combinatorial Calabi energy}\ls

Similar to \cite{Ge1}, the combinatorial Calabi energy in hyperbolic background geometry is defined as
\begin{equation}
\mathcal{C}(r)=\|K\|^2=\sum_{i=1}^{N}K_{i}^2.
\end{equation}

\begin{rmk}
In Euclidean background geometry, the average curvature $K_{av}=\sum_{i=1}^{N}K_{i}/N=2\pi \chi(X)/N$ is a combinatorial invariant. Then
\begin{equation*}
\sum_{i=1}^{N}(K_{i}-K_{av})^2=\sum_{i=1}^{N}K_{i}^2-NK_{av}^2=\sum_{i=1}^{N}K_{i}^2-\frac{(2\pi \chi(X))^2}{N}.
\end{equation*}
Therefore, it makes no difference to define combinatorial Calabi energy as $\mathcal{C}(r)=\sum_{i=1}^{N}K_{i}^2$ or as $\mathcal{C}(r)=\sum_{i=1}^{N}(K_{i}-K_{av})^2$, since they differ only by a combinatorial constant.
However, in hyperbolic background geometry, the average curvature $K_{av}=\frac{2\pi \chi(X)+Area(X)}{N}$
contains the area term and hence is no longer a combinatorial invariant.
Therefore, the combinatorial Calabi energy we define here do not contain $K_{av}$ term.
\end{rmk}

Set $u_i=\ln \tanh \frac{r_i}{2},$ then the map $u=u(r)$ maps $\mathds{R}^{N}_{>0}$ to $\mathds{R}^{N}_{<0}$ homeomorphically. In the following we sometimes treat functions in $r_i$, and sometimes in $u_i$.

The gradient of the combinatorial Calabi energy $\mathcal{C}(u)=\sum_{i=1}^{N}K_i^2$ is
\begin{equation*}
 \begin{aligned}
 \nabla_u \mathcal{C}
 =\begin{pmatrix}
                      \nabla_{u_1}\mathcal{C} \\
                      \vdots \\
                      \nabla_{u_N}\mathcal{C} \\
                    \end{pmatrix}
 =2\begin{pmatrix}
     \frac{\partial K_1}{\partial u_1} & \cdots & \frac{\partial K_N}{\partial u_1} \\
     \vdots & \ddots & \vdots \\
     \frac{\partial K_1}{\partial u_N} & \cdots & \frac{\partial K_N}{\partial u_N} \\
   \end{pmatrix}
   \begin{pmatrix}
     K_1 \\
     \vdots \\
     K_N \\
   \end{pmatrix}
 =2L^TK,
 \end{aligned}
\end{equation*}
where $K=(K_1, \cdots, K_N)^T$ and
\begin{eqnarray}\label{hyper-discrete-Laplacian}
  L=(L_{ij})_{N\times N}=\frac{\partial (K_1,\cdots,K_N)}{\partial (u_1,\cdots,u_N)}
  =\begin{pmatrix}
     \frac{\partial K_1}{\partial u_1} & \cdots & \frac{\partial K_1}{\partial u_N} \\
     \vdots & \ddots & \vdots \\
     \frac{\partial K_N}{\partial u_1} & \cdots & \frac{\partial K_N}{\partial u_N} \\
   \end{pmatrix}
\end{eqnarray}
is the Jacobian of the curvature map $K=K(u)$. In \cite{Ge1} and \cite{Ge2}, we had interpret the Jacobian of the curvature map $K$ to be a type of discrete Laplace operator. Here, (\ref{hyper-discrete-Laplacian}) is also a type of discrete Laplacian, called  discrete dual-Laplacian in hyperbolic background geometry and written as $\mathds{H}^2$-DDL for short.

\subsection{The Combinatorial Calabi flow}
\begin{defn}
 Given $(X, T, \Phi,\mathds{H}^2)$, where $X$, $T$, $\Phi\in[0, \pi/2]$ are defined as before.
 The combinatorial Calabi flow in hyperbolic geometry background is defined as
\begin{equation}\label{Calabi-flow-hyper}
\dot{u}(t)=-L^TK.
\end{equation}
\end{defn}

\begin{prop}
Combinatorial Calabi flow $(\ref{Calabi-flow-hyper})$ is the negative gradient flow of Combinatorial Calabi energy and the solution of $(\ref{Calabi-flow-hyper})$ exists locally. Moreover, combinatorial Calabi energy is descending along the flow.
\end{prop}
\pf By the definition of combinatorial Calabi flow (\ref{Calabi-flow-hyper}), we have $\dot{u}(t)=-L^TK=-\frac{1}{2}\nabla_u \mathcal{C}$.
The flow is a system of ODEs. Therefore, given any initial metric $u(0)\in \mathds{R}^{N}_{<0}$, or equivalently $r(0)\in \mathds{R}^{N}_{>0}$,
the solution of Calabi flow exists for the interval $t\in [0,\epsilon)$, where $\epsilon$ is a positive number small enough.
By direct computation, we have $\mathcal{C}'(t)=-\frac{1}{2}\|\nabla_u \mathcal{C}\|^2\leq 0$.
Thus the combinatorial Calabi energy is descending along the flow. \qed

\begin{prop}
 The combinatorial curvature evolves according to $$\dot{K}=-LL^TK.$$
\end{prop}
\pf $\dot{K}=\frac{\partial (K_1,\cdots,K_N)}{\partial (u_1,\cdots,u_N)}\dot{u}=L\dot{u}=-LL^TK.$ \qed

Denote $\Gamma(u)=-L^TK$, then the combinatorial Calabi flow equation (\ref{Calabi-flow-hyper}) could be written as
$\dot{u}=\Gamma(u)$,
which is an autonomous system. The function $\Gamma(u)$ (or the matrix $L$) determines all the properties of the autonomous system.

\section{Discrete Dual-Laplacian}\ls
\subsection{Discrete dual-Laplacian}\ls
We derive the explicit form of $L_{ij}=\frac{\partial K_i}{\partial u_j}$ first.
As $\frac{\partial \theta_{i}^{jk}}{\partial u_j}=\frac{\partial \theta_{j}^{ik}}{\partial u_i}$
and $\frac{\partial K_i}{\partial u_j}=\frac{\partial K_j}{\partial u_i}$ \cite{CL1},
we have $L_{ij}=L_{ji}$. We write $j\thicksim i$ in the following if the vertices $i$ and $j$ are adjacent.
\begin{enumerate}
  \item If $j\thicksim i$, we have
  \begin{equation*}
    \begin{aligned}
        \frac{\partial K_i}{\partial u_j}
      =&\frac{\partial K_i}{\partial r_j}\sinh r_j\\
      =&\frac{\partial\left(2\pi-\sum_{\{i,j,k\}\in F}\theta_{i}^{jk}\right)}{\partial r_j}\sinh r_j\\
      =&-\sum_{\{i,j,k\}\in F}\frac{\partial \theta_{i}^{jk}}{\partial r_j}\sinh r_j\\
      =&-\frac{\partial (\theta_{i}^{jk}+\theta_{i}^{jl})}{\partial r_j}\sinh r_j,
    \end{aligned}
  \end{equation*}
  where $k,l$ are the vertices such that $\{i,j,k\}$,$\{i,j,l\}$ are adjacent faces.
 Set $B_{ij}=\frac{\partial \theta_{i}^{jk}}{\partial r_j}\sinh r_j+\frac{\partial \theta_{i}^{jl}}{\partial r_j}\sinh r_j$, then
  \begin{equation}\label{jiadjacent}
    \frac{\partial K_i}{\partial u_j}=-B_{ij}\,,\,\,j\thicksim i.
  \end{equation}
  \item If $j=i$, we have
  \begin{equation*}
    \begin{aligned}
        \frac{\partial K_i}{\partial u_i}
      =&\frac{\partial K_i}{\partial r_i}\sinh r_i\\
      =&\frac{\partial\left(2\pi-\sum_{\{i,j,k\}\in F}\theta_{i}^{jk}\right)}{\partial r_i}\sinh r_i\\
      =&-\sum_{\{i,j,k\}\in F}\frac{\partial \theta_{i}^{jk}}{\partial r_i}\sinh r_i\\
      =&\sum_{\{i,j,k\}\in F}\frac{\partial \left(\theta_{j}^{ik}+\theta_{k}^{ij}+Area(\triangle v_iv_jv_k)\right)}{\partial r_i}\sinh r_i\\
      =&\sum_{\{i,j,k\}\in F}\left(\frac{\partial \theta_{i}^{jk}}{\partial r_j}\sinh r_j+\frac{\partial \theta_{i}^{kj}}{\partial r_k}\sinh r_k\right)+\sum_{\{i,j,k\}\in F}\frac{\partial Area(\triangle v_iv_jv_k)}{\partial r_i}\sinh r_i\\
      =&\sum_{j\thicksim i}\left(\frac{\partial \theta_{i}^{jk}}{\partial r_j}\sinh r_j+\frac{\partial \theta_{i}^{jl}}{\partial r_j}\sinh r_j\right)+\sinh r_i\frac{\partial}{\partial r_i}\left(\sum_{\{i,j,k\}\in F} Area(\triangle v_iv_jv_k)\right).
    \end{aligned}
  \end{equation*}
  Set
    $$A_i=\sinh r_i\frac{\partial}{\partial r_i}\left(\sum_{\{i,j,k\}\in F} Area(\triangle v_iv_jv_k)\right),$$
then we have
    \begin{equation}\label{jisame}
      \frac{\partial K_i}{\partial u_i}= A_i+\sum_{j\thicksim i}B_{ij}.
    \end{equation}
  \item If $j\nsim i$ and $j \neq i$, then
    \begin{equation}\label{jielse}
      \frac{\partial K_i}{\partial u_j}=0.
    \end{equation}
\end{enumerate}

Set $A=diag\big\{A_1,\cdots ,A_N\big\}$ and $L_B=\big((L_B)_{ij}\big)_{N\times N}$, where
\begin{gather*}
\left(L_B\right)_{ij}=
\begin{cases}
\,\,\sum\limits_{k \sim i}B_{ik} \,, & \text{$ j=i,$} \\
\,\,\,\,\,-B_{ij} \,,& \text{$ j\sim i,$} \\
\,\,\,\,\,\,\,\,\,0 \,,& \text{$ j\nsim i,\, j\neq i,$}
\end{cases}
\end{gather*}
then we have
\begin{thm}\label{positivity of L}
 The matrix $L$ is positive definite and
$$L=A+L_B,$$
where $A$ is positive definite and $L_B$ is semi-positive definite.
\end{thm}
\pf
Combining the formulas (\ref{jiadjacent}), (\ref{jisame}) and (\ref{jielse}), we get
\begin{gather*}
L_{ij}=
\begin{cases}
A_i+\sum\limits_{k \sim i}B_{ik} \,, & \text{$ j=i,$} \\
\,\,\,\,\,\,\,-B_{ij} \,,& \text{$ j\sim i,$} \\
\,\,\,\,\,\,\,\,\,\,\,\,0 \,,& \text{$ j\nsim i,\, j\neq i.$}
\end{cases}
\end{gather*}
By Lemma 2.2. in \cite{CL1}, we know $A_i>0$, $B_{ij}>0$ for $j\sim i$.
By Lemma 3.10. in \cite{CL1}, we know the matrix $L_B$ is semi-positive definite.
Thus $L$ is positive definite. \qed

The classical discrete Laplace operator $``\Delta"$ is often written as the following form (\cite{CHU}) 
\begin{equation*}
\Delta f_i=\sum_{j\sim i}\omega_{ij}(f_j-f_i).
\end{equation*}
Notice that the weight $\omega_{ij}$ can be arbitrarily selected for different purpose. Using a special weight which comes from dual structure of circle patterns, we established $\Delta=-L$ in \cite{Ge1} and \cite{Ge2}. Formally, we also define discrete Laplace operator $\Delta=-L$ in this paper. Both $\Delta$ and $L$ act on function $f$ by matrix multiplication, that is
\begin{equation*}
\Delta f=-Lf,
\end{equation*}
or in component form
\begin{equation*}
\Delta f_i=\sum\limits_{j \sim i}B_{ij}(f_j-f_i)-A_if_i.
\end{equation*}
It's interesting that the combinatorial Calabi flow (\ref{Calabi-flow-hyper}) can be written as
\begin{equation*}
\cfrac{dr_i}{dt}=\sinh r_i \Delta K_i\ \ \ \text{or}\ \ \ \ \cfrac{du_i}{dt}=\Delta K_i.
\end{equation*}
Its matrix form $\cfrac{du}{dt}=\Delta K$ is similar to the 2-dimension smooth Calabi flow $\cfrac{\partial g}{\partial t}=\Delta K$.

\begin{rmk}\label{global rigidity}
(Global Rigidity) Consider the combinatorial Ricci potential
\begin{equation*}
f(u)\triangleq \int_{u_0}^u\sum_{i=1}^NK_idu_i\,,\,\,u\in \mathds{R}^N_{<0},
\end{equation*}
where $u_0$ is an arbitrarily point in $\mathds{R}^N_{<0}$. The integral is well-defined,
since $\sum_{i=1}^NK_idu_i$ is a closed differential form.
Moreover, because $Hess(f)=\frac{\partial(K_1,...,K_N)}{\partial(u_1,...,u_N)}=L$ is positive definite by Theorem $\ref{positivity of L}$,
$\nabla f=K:\mathds{R}^N_{<0}\rightarrow\mathds{R}^N$ is an embedding. Thus the combinatorial curvature determines $u$ and hence the circle packing metric $r$ uniquely. This property is called global rigidity  \cite{CL1, T1}.
\end{rmk}

\subsection{Some estimates of the discrete Laplacian with $\Phi\equiv0$}
We will establish some estimates of the elements of discrete Laplacian if the weight function $\Phi\equiv0$.
\begin{lem} \label{estimate-lemma}
For any $x, y, z > 0$,
\begin{equation}
0<\frac{\sinh(2x+y+z)}{\sinh(x+y)\sinh(x+z)}\sqrt{\frac{\sinh x \sinh y \sinh z }{\sinh(x+y+z)}}<\cosh 1.
\end{equation}
\end{lem}
As the proof of the  lemma is elementary and tedious, we postpone it to Appendix \ref{appendix-A}.
\begin{coro} \label{estimate-corollary}
Suppose a hyperbolic triangle $\triangle ijk$ is configured by circle packing metric $r_i, r_j, r_k > 0$ with zero weight. Then
\begin{equation*}
\begin{aligned}
&\ \ \ 0<\frac{\partial \theta^{jk}_i}{\partial r_j}\sinh r_j<\frac{1}{2},\\
&0<-\frac{\partial \theta^{jk}_i}{\partial r_i}\sinh r_i<\cosh 1,\\
0<&\frac{\partial Area(\triangle ijk)}{\partial r_i}\sinh r_i<\cosh 1.
\end{aligned}
\end{equation*}
\end{coro}
\pf If $\Phi\equiv0$, we can derive the following expressions directly by a lengthy but ordinary calculation using hyperbolic cosine law:
\begin{equation*}
\frac{\partial \theta_{i}^{jk}}{\partial r_j}\sinh r_j
=\cfrac{1}{\sinh (r_i+r_j)}\sqrt{\frac{\sinh r_i\sinh r_j\sinh r_k}{\sinh (r_i+r_j+r_k)}},
\end{equation*}
\begin{equation*}
\begin{aligned}
\frac{\partial \theta_{i}^{jk}}{\partial r_i}\sinh r_i
=-\frac{\sinh(2r_i+r_j+r_k)}{\sinh(r_i+r_j)\sinh(r_i+r_k)}\sqrt{\frac{\sinh r_i\sinh r_j\sinh r_k}{\sinh (r_i+r_j+r_k)}}.
\end{aligned}
\end{equation*}
Therefore, we have
$$0<\frac{\partial \theta^{jk}_i}{\partial r_j}\sinh r_j<\frac{1}{2}.$$
By Lemma \ref{estimate-lemma}, we have
$$0<-\frac{\partial \theta^{jk}_i}{\partial r_i}\sinh r_i<\cosh 1.$$
By the formula $\theta^{jk}_i+\theta^{ik}_j+\theta^{ij}_j=\pi-Area(\triangle ijk)$ and Lemma \ref{estimate-lemma},
we have
\begin{equation*}
\begin{split}
0 & <\frac{\partial Area(\triangle ijk)}{\partial r_i}\sinh r_i \\
  & = \frac{\sinh(2r_i+r_j+r_k)-\sinh(r_i+r_j)-\sinh(r_i+r_k)}{\sinh(r_i+r_j)\sinh(r_i+r_k)}
      \sqrt{\frac{\sinh r_i \sinh r_j \sinh r_k }{\sinh(r_i+r_j+r_k)}} \\
  & < \cosh 1.
\end{split}
\end{equation*}
\qed
\begin{rmk}
 For arbitrary weight $\Phi: E\rightarrow[0,\pi/2]$, Bennett Chow and Luo Feng \cite{CL1} had got that
\begin{equation*}
\begin{aligned}
\frac{\partial \theta_{i}^{jk}}{\partial r_j}\sinh r_j
=&\frac{\sinh r_j}{\sinh l_{ij}}\frac{\cos \theta_{j}^{ki'}-\cos \theta_{j}^{ik'}\cos \theta_{j}^{ik}}
  {\sin \theta_{j}^{ik}},
\end{aligned}
\end{equation*}
\begin{equation*}
\begin{aligned}
\frac{\partial \theta_{i}^{jk}}{\partial r_i}\sinh r_i
=-\frac{\sinh r_i}{\sinh l_{ij}}\frac{\cos\theta_{j}^{ik}\cos \theta_{i}^{jk'}+\cos\theta_{k}^{ij}\cos \theta_{i}^{kj'}}
 {\sin \theta_{j}^{ik}}.
\end{aligned}
\end{equation*}
\end{rmk}

By Corollary \ref{estimate-corollary}, we can get the following proposition easily.
\begin{prop} \label{basic-estimate}
For fixed $(X,T,\Phi\equiv0,\mathds{H}^2)$, we have
\begin{equation*}
0<B_{ij}<1,
\end{equation*}
\begin{equation*}
0<A_i<d_i\cosh 1\leq d\cosh 1,
\end{equation*}
where $d=\max\limits_{1\leq i\leq N}\left\{d_i\right\}$, $d_i$ is the degree at vertex $i$.
\end{prop}

\begin{coro} \label{r-lower-bound-0-weight}
Given $(X, T, \Phi\equiv0, \mathds{H}^2)$, then $r(t)$, the solution  of the combinatorial Calabi flow $(\ref{Calabi-flow-hyper})$,
has a positive lower bound in finite time.
\end{coro}
\pf
Note that $(2-d)\pi<K_i<2\pi$ and $\Delta K_i=\sum\limits_{j \sim i}B_{ij}(K_j-K_i)-A_iK_i$,
then by combinatorial Calabi flow equation
$\frac{du_i}{dt}=\Delta K_i$
and the estimates of $B_{ij}$ and $A_i$ in Proposition \ref{basic-estimate}, we have
$$|u_i(t)-u_i(0)|<c_2t,$$
where $c_2$ is a constant depending on the combinatorial structure of the triangulation.
As $u_i=\ln \tanh \frac{r_i}{2}$, we have
\begin{equation}\label{estimate of r_i(t) with weight 0}
0<c_1e^{-c_2t}<\tanh \frac{r_i(t)}{2}<c_3e^{c_2t},
\end{equation}
where $c_1=\min\limits_{1\leq i\leq N}\left\{\tanh \frac{r_i(0)}{2}\right\}$,
$c_3=\max\limits_{1\leq i\leq N}\left\{\tanh \frac{r_i(0)}{2}\right\}$.
Notice that $c_3e^{c_2t}\rightarrow+\infty$ as $t\rightarrow+\infty$, while $\tanh$ is always upper bounded by $1$,
the estimate on the right hand side of $(\ref{estimate of r_i(t) with weight 0})$ will be of no use if $t$ is big enough.
Thus we can only get an estimate of the lower bound of $r_i(t)$,
\begin{equation*}
r_i(t)\geq \ln \cfrac{1+c_1e^{-c_2t}}{1-c_1e^{-c_2t}}>0.
\end{equation*}
The estimate shows that the solution $r(t)$ of the combinatorial Calabi flow (\ref{Calabi-flow-hyper})
has a positive lower bound in finite time. Specially, if $t\in [0,T]$, then we have $r_i(t)>c_T>0$, where $c_T=\ln \cfrac{1+c_1e^{-c_2T}}{1-c_1e^{-c_2T}}>0$.\qed

\section{Local Convergence of Combinatorial Calabi Flow}
\subsection{ZCCP-metric (zero curvature circle packing metric)}

Zero curvature circle packing metric is a circle packing metric that determines constant zero curvature, written as ZCCP-metric for short.
we denote it as $r^*$ (or $u^*$ in $u$-coordinate) if it exists. The ZCCP-metric is important
for the combinatorial Calabi flow (\ref{Calabi-flow-hyper}),
because it is the unique stable point of the autonomous system $\dot{u}=\Gamma(u)=-LK$.

There are combinatorial obstructions for the existence of ZCCP-metric, which was found by Thurston \cite{T1}. We quote the combinatorial conditions here.
\begin{prop}
For fixed $(X,T,\Phi,\mathds{H}^2)$, where $X$, $T$, $\Phi\in[0, \pi/2]$ are defined as before,
there exists ZCCP-metric if and only if the following two combinatorial conditions are satisfied simultaneously:
\begin{enumerate}
  \item For any three edges $e_1, e_2, e_3$ forming a null homotopic loop in $X$, if $\sum_{i=1}^3\Phi(e_i)\geq\pi$, then $e_1, e_2, e_3$ form the boundary of a triangle of $T$;
  \item For any four edges $e_1, e_2, e_3, e_4$ forming a null homotopic loop in $X$, if $\sum_{i=1}^4\Phi(e_i)\geq2\pi$, then $e_1, e_2, e_3, e_4$ form the boundary of the union of two adjacent triangles.
\end{enumerate}
\end{prop}

\subsection{Convergence with small initial Calabi energy }
\begin{lem} \label{compact-converge}
Given $(X,T,\Phi,\mathds{H}^2)$, where $X$, $T$, $\Phi\in[0, \pi/2]$ are defined as before.
If the solution of the combinatorial Calabi flow $(\ref{Calabi-flow-hyper})$ lies in a compact subset of $\mathds{R}^{N}_{<0}$,
then the solution exists for all $t\in [0,+\infty)$ and the solution has exponential convergence rate.
\end{lem}
\pf Denote the solution by $u(t)$. $u(t)$ must exists for all $t\in[0,+\infty)$, otherwise $u(t)$ will attain the boundary of $\mathds{R}^{N}_{<0}$.
This contradicts the condition $\{u(t)\}\subset\subset\mathds{R}^{N}_{<0}$.
 By the same condition, the eigenvalue of matrix $L$ has a uniform positive lower bound $\lambda_1$,
 that is, $\lambda\big(L(t)\big)\geq \lambda_1$ along the combinatorial Calabi flow. Hence
$$\mathcal{C}'(t)=2K^T \dot{K}(t)=-2K^TL^2K\leq -2\lambda^2_1K^TK=-2\lambda^2_1\mathcal{C}(t).$$
So $\mathcal{C}(t)\leq\mathcal{C}(0)\exp^{-2\lambda^2_1t}$ and $\big|K_i(t)\big| \leq \big|K(t)\big|=\sqrt{\mathcal{C}(t)}\leq \sqrt{\mathcal{C}(0)}\exp^{-\lambda^2_1t}$. As $\{u(t)\}\subset\subset\mathds{R}^{N}_{<0}$,
we know $L_{ij}$ and $\sinh r_i$ are bounded along the combinatorial Calabi flow. Hence
\begin{equation*}
\begin{aligned}
\Big|\frac{du_i}{dt}\Big|&=\big|\triangle K_i\big|=\Big|\sum_jL_{ij}K_j\Big|\leq c_1\exp^{-\lambda^2_1t},\\
\Big|\frac{dr_i}{dt}\Big|&=\big|\sinh r_i\triangle K_i\big|=\Big|\sinh r_i \sum_jL_{ij}K_j\Big|\leq c_2\exp^{-\lambda^2_1t},
\end{aligned}
\end{equation*}
where $c_1$, $c_2$ are positive constants.
This implies that the solution converges with exponential rate.
\qed

Assuming the Calabi flow exists for all time and converges, then we have
\begin{prop} \label{converg-zero-metric}
Given $(X,T,\Phi,\mathds{H}^2)$, where $X$, $T$, $\Phi\in[0, \pi/2]$ are defined as before.
If the solution of Calabi flow $(\ref{Calabi-flow-hyper})$ exists for $t\in [0,+\infty)$
and converges for some initial circle packing metric, then
\begin{enumerate}
  \item the solution has exponential convergence rate;
  \item there exists ZCCP-metric on $(X,T,\Phi,\mathds{H}^2)$;
  \item the Euler characteristic of $X$ is negative, \emph{i.e.} $\chi(X)<0$.
\end{enumerate}
\end{prop}
\pf
As the solution of combinatorial Calabi flow converges for some initial data, the solution $u(t)$ will lie
in a compact subset in $\mathds{R}^{N}_{<0}$ from some time $T$. Thus the solution has
exponentially convergence rate by Lemma \ref{compact-converge}.

Suppose $r(t), t\in [0,+\infty)$, is the solution to the combinatorial flow,
$K(t)$ and $L(t)$ are the corresponding curvature and discrete Laplacian along the flow.
By the hypothesis in the proposition, we have $r(+\infty)=\lim_{t\rightarrow +\infty}r(t)\in \mathds{R}^{N}_{>0}$
exists, so does $K(+\infty)=\lim_{t\rightarrow +\infty}K(t)$ and $L(+\infty)=\lim_{t\rightarrow +\infty}L(t)$.
Thus the combinatorial Calabi energy $\mathcal{C}(t)=\sum_{i=1}^{N}K_{i}^{2}(t)$ converges as $t\rightarrow\infty$.
The derivative of the combinatorial Calabi energy
$$\mathcal{C}'(t)=2K^T \dot{K}(t)=-2K^TL^2K$$
converges also. As $\mathcal{C}'(t)\leq 0$, we have $\mathcal{C}'(+\infty)=\lim_{t\rightarrow +\infty}{C}'(t)\leq 0$.
Then $\mathcal{C}(+\infty)$ exists with $\mathcal{C}(+\infty)\geq 0$ and
$\mathcal{C}'(+\infty)$ exists with $\mathcal{C}'(+\infty)\leq 0$.
This implies that $\mathcal{C}'(+\infty)=0$, which is equivalent to
$$K^TL^2K(+\infty)=0.$$
By the positivity of $L$, we have $K(+\infty)=0$. Thus $r(+\infty)$ is a ZCCP-metric.
By the combinational Gauss-Bonnet formula (\ref{Gauss-Bonnet formula}), we have $\chi(X)<0$.
\qed

\begin{thm}
For fixed $(X,T,\Phi,\mathds{H}^2)$, where $X$, $T$, $\Phi\in[0, \pi/2]$ are defined as before,
suppose there exists a circle packing metric $u^*$ that determines constant zero curvature,
then the solution of combinatorial Calabi flow $(\ref{Calabi-flow-hyper})$ exists for all $t\in [0,+\infty)$ and
converges exponentially fast to $u^*$
if the initial combinatorial Calabi energy is small enough.
\end{thm}
\pf We claim that $u^*$ is the only asymptotically stable point of the combinatorial Calabi flow (\ref{Calabi-flow-hyper}).
 In fact,
$$D_{u^*}(-LK)=-L^2(u^*)<0.$$
So the claim is true by the Lyapunov stability theorem. This implies the solution of the combinatorial Calabi flow (\ref{Calabi-flow-hyper})
exists for $t\in [0,+\infty)$ and converges exponentially to $u^*$ if the initial metric $u(0)$ is close to $u^*$.
By Remark \ref{global rigidity}, this is equivalent to the initial combinatorial Calabi energy is small enough.\qed

Similarly, we can prove
\begin{thm}
Given any user prescribed combinatorial curvature $\bar{K}$ on fixed $(X,T,\Phi,\mathds{H}^2)$,
where $X$, $T$, $\Phi\in[0, \pi/2]$ are defined as before.
If $\bar{K}$ is admissible, that is, there exists a circle packing metric $\bar{r}$ with $\bar{K}=K(\bar{r})$,
and the initial modified combinatorial Calabi energy $\|K(0)-\bar{K}\|^2$ is small enough,
then the solution of the following modified combinatorial Calabi flow
\begin{equation}\label{modified Calabi flow}
\dot{u}=L(\bar{K}-K)
\end{equation}
exists for $t\in [0,+\infty)$ and
converges exponentially fast to $\bar{u}$.
\end{thm}
The modified combinatorial Calabi flow (\ref{modified Calabi flow}) provides an algorithm to compute circle packing metrics with user prescribed combinatorial curvatures. In fact, any algorithm, aiming at minimizing the modified combinatorial Calabi energy $\|K-\bar{K}\|^2$, tends to find circle packing metric $\bar{r}$ automatically.

\section{Global convergence for $K_i<C<2\pi$}
Along the combinatorial Calabi flow, the evolution equation of combinatorial curvature $\dot{K}=-\Delta^2K$ contains the $\Delta^2$ term. It is a fourth order equation. There are not efficient ``Maximal Principle" for fourth order equations, this may be the most important differences between the combinatorial Calabi flow and the combinatorial Ricci flow introduced in \cite{CL1}. We give some global convergence results assuming that the curvatures are uniformly bounded from above by a constant smaller than $2\pi$.

\begin{lem} \label{estimate-Ricci-potential}
For fixed $(X, T, \Phi, \mathds{H}^2)$, where $X$, $T$, $\Phi\in[0, \pi/2]$ are defined as before.
The combinatorial Ricci potential is defined as
\begin{equation*}
f(u)= \int_{u_0}^u\sum_{i=1}^N\big(K_i-K_i(u_0)\big)du_i\,,\,\,u\in \mathds{R}^N_{<0}\,,
\end{equation*}
where $u_0$ is an arbitrarily point in $\mathds{R}^N_{<0}$. Then
$\lim\limits_{\|u\|\rightarrow +\infty}f(u)=+\infty.$
\end{lem}
Since the proof is so tedious, we defer it to the Appendix \ref{appendix B}. Using Lemma \ref{estimate-Ricci-potential}, we get
\begin{thm} \label{r-lower-bound}
Given $(X, T, \Phi, \mathds{H}^2)$, where $X$, $T$, $\Phi\in[0, \pi/2]$ are defined as before.
Then $r(t)$, the solution  of Calabi flow $(\ref{Calabi-flow-hyper})$,
has a positive lower bound in finite time.
\end{thm}
\pf
Let $u(t)$ be the solution of $(\ref{Calabi-flow-hyper})$ in $u$-coordinate. Denote $\varphi(t)=f(u(t))$, then $\varphi'(t)=-K^TLK\leq0$, $\varphi$ is descending and  $\varphi(t)\leq \varphi(0)$, hence $\{u(t)\}\subset f^{-1}[0,\varphi(0)]$. Using Lemma \ref{estimate-Ricci-potential}, we know $|u_i(t)|$ is uniformly bounded from above. Moreover, using similar estimate in the proof of Corollary \ref{r-lower-bound-0-weight}, we get the lower bound of $r_i(t)$.
\qed

\begin{thm} \label{0-weight-final-them}
For $(X, T, \Phi, \mathds{H}^2)$, where $X$, $T$, $\Phi\in[0, \pi/2]$ are defined as before.
Assuming the curvature $K_i$ at each vertex $v_i$ has a uniform upper bound smaller than $2\pi$ along the combinatorial Calabi flow $(\ref{Calabi-flow-hyper})$, then
\begin{enumerate}
  \item the solution $r(t)$ of $(\ref{Calabi-flow-hyper})$ exists for all $t\in[0,+\infty)$;
  \item $r(t)$ converges if and only if ZCCP-metric $r^*$ exists;
  \item If $r(t)$ converges, then it converges to $r^*$ exponentially and $\chi(X)<0$.
\end{enumerate}
\end{thm}
\pf
Suppose $r(t)$ is the solution of (\ref{Calabi-flow-hyper}) on the maximal time interval $t\in[0,T)$. If $T<+\infty$, the metric will blow up near $T$, that is, there are times $t_n\uparrow T$ such that $r_i(t_n)\rightarrow +\infty$ at some vertex $v_i$, or there are times $t_n\uparrow T$ such that $r_j(t_n)\rightarrow 0$ at some vertex $v_j$. We show that these two ``blow up" phenomena will never occur. Using hyperbolic cosine law in triangles which has $v_i$ as one of its vertices, we can get $\theta_i^{jk}(t_n)\rightarrow 0$, thus $K_i(t_n)\rightarrow 2\pi$. This contradicts the hypothesis. Using Corollary \ref{r-lower-bound-0-weight}, all $r_i$ must bounded from below and will never tend to zero in finite time $t\leq T$. Hence $r(t)$ exists for all $t\in[0,+\infty)$.

Suppose constant zero curvature metric $r^*$ exists. we want to prove $u(t)\subset\subset\mathds{R}^N_{<0}$.
Let us consider the combinatorial Ricci potential
\begin{equation*}
f(u)= \int_{u^*}^u\sum_{i=1}^NK_idu_i\,,\,\,u\in \mathds{R}^N_{<0},
\end{equation*}
where $u^*$ is zero curvature metric in $u$-coordinate. Using Theorem \ref{r-lower-bound}, we know $|u(t)|$ is bounded above. If $u(t)$ is not contained in a relative compact set of $\mathds{R}^N_{<0}$, then there are times $t_n\uparrow +\infty$ such that $u_i(t_n)\rightarrow0$ at some vertex $v_i$, hence $K_i(t_n)\rightarrow 2\pi$. This contradicts the hypothesis. Therefore, $u(t)\subset\subset\mathds{R}^N_{<0}$.

The proof can be finished by using
 Lemma \ref{compact-converge} and Lemma \ref{converg-zero-metric}. \qed

\begin{rmk}
Combinatorial Calabi flow in $\mathbb{S}^2$-geometry.

In spherical background geometry case, all triangles are embedded in the standard $\mathbb{S}^2$. The admissible circle packing metrics compose the space
$$\mathcal {M}_r\doteqdot \left\{r=(r_1,...,r_N)^T\in \mathds{R}^N_{>0} \;\Big|\;r_i+r_j+r_k<\pi,\,\forall\triangle v_iv_jv_k\in F \right\}.$$
The coordinate transformations are $u_i=\ln \tan \frac{r_i}{2}$, and the corresponding $u$-metric space is
$\mathcal {M}_u=\left\{u\in \mathds{R}^N \Big|u_i=\ln \tan \frac{r_i}{2}, \;r\in \mathcal {M}_r \right\}$. $\mathcal {M}_r$ is open and simply connected, but not convex.
The Jacobian of curvature map $K=K(u)$ is still denoted by $L$. Then $L=L_B-A$, where $L_B$ and $A$ are defined similarly (using $\sin$ instead of $\sinh$). Also, $L_B$ is semi-positive definite and $A$ is positive definite.
The combinatorial Calabi flow in $\mathbb{S}^2$-geometry is defined to be $$\frac{du}{dt}=\Delta K=-LK=-\frac{1}{2}\nabla_u \mathcal{C}$$ in matrix form, or $$\cfrac{dr_i}{dt}=\sin r_i \Delta K_i$$ in component form. The solution of this flow exists locally and the Calabi energy is descending along this flow. However, the signature of $L$, the number of stable points (that is, the solution of the matrix equation $\Gamma(u)=(A-L_B)K=0$) besides $K=0$,
and the convergence behavior of this flow are unclear.
\qed
\end{rmk}

\appendix
\section{A Proof of Lemma \ref{estimate-lemma}}\label{appendix-A}
\begin{lem} \label{appendix-lemma1}
For any $x, y, z > 0$
\begin{equation*}
0<\sqrt{\frac{\sinh x \sinh y \sinh z }{\sinh(x+y+z)}}<\frac{1}{2}.
\end{equation*}
\end{lem}
\pf
For any $x,y>0$
$$0< 2\sinh x \sinh y < \sinh x \cosh y + \cosh x \sinh y = \sinh(x+y).$$
Thus for any $x,y,z>0$, we have
$$0< 4 \sinh x \sinh y \sinh z <2\sinh x \sinh(y+z)<\sinh(x+y+z),$$
that is
$$0<\frac{\sinh x \sinh y \sinh z }{\sinh(x+y+z)}<\frac{1}{4}.$$
\qed

\begin{lem} \label{appendix-lemma2}
For any $x, y, z > 0$
\begin{equation*}
0<\frac{1}{\sinh(x+y)}\sqrt{\frac{\sinh x \sinh y \sinh z }{\sinh(x+y+z)}}<\frac{1}{2}.
\end{equation*}
\end{lem}
\pf
As
$\sinh(x+y)\geq \sinh x +\sinh y , \forall x,y>0$
and $0<\sinh z< \sinh(x+y+z)$, we have
$$0<\frac{1}{\sinh(x+y)}\sqrt{\frac{\sinh x \sinh y \sinh z }{\sinh(x+y+z)}}<\frac{\sqrt{\sinh x \sinh y}}{\sinh(x+y)}\leq
\frac{1}{2}\frac{\sinh x+ \sinh y}{\sinh(x+y)}\leq \frac{1}{2}.$$
\qed

\begin{lem} \label{appendix-lemma3}
For any $x, y, z > 0$
\begin{equation*}
0<\frac{\cosh(x+y)}{\sinh(x+y)}\sqrt{\frac{\sinh x \sinh y \sinh z }{\sinh(x+y+z)}}<\frac{\cosh1}{2}.
\end{equation*}
\end{lem}
\pf
For any $x, y, z > 0$,
if $x+y\geq1$, using Lemma \ref{appendix-lemma1}, we get
$$0<\frac{\cosh(x+y)}{\sinh(x+y)}\sqrt{\frac{\sinh x \sinh y \sinh z }{\sinh(x+y+z)}}<\frac{1}{2}\frac{1}{\tanh 1}<\frac{\cosh1}{2},$$
since function $\tanh x$ is strictly increasing on $(0,+\infty)$;
if $x+y<1$, using Lemma \ref{appendix-lemma2}, we get
$$0<\frac{\cosh(x+y)}{\sinh(x+y)}\sqrt{\frac{\sinh x \sinh y \sinh z }{\sinh(x+y+z)}}<\frac{\cosh1}{2},$$
since function $\cosh x$ is strictly increasing on $(0,+\infty)$.
\qed

\begin{lem}
For any $x, y, z > 0$,
\begin{equation*}
0<\frac{\sinh(2x+y+z)}{\sinh(x+y)\sinh(x+z)}\sqrt{\frac{\sinh x \sinh y \sinh z }{\sinh(x+y+z)}}<\cosh 1.
\end{equation*}
\end{lem}
\pf
As
$$\frac{\sinh(2x+y+z)}{\sinh(x+y)\sinh(x+z)}=\frac{\cosh(x+y)}{\sinh(x+y)}+\frac{\cosh(x+z)}{\sinh(x+z)},$$
the above estimate follows by  Lemma \ref{appendix-lemma3}.
\qed

\section{B Proof of Lemma \ref{estimate-Ricci-potential}}\label{appendix B}
\begin{lem}
For fixed $(X, T, \Phi, \mathds{H}^2)$, where $X$, $T$, $\Phi\in[0, \pi/2]$ are defined as before.
The combinatorial Ricci potential is defined as
\begin{equation*}
f(u)\triangleq\int_{u_0}^u\sum_{i=1}^N\big(K_i-K_i(u_0)\big)du_i\,,\,\,u\in \mathds{R}^N_{<0}\,,
\end{equation*}
where $u_0$ is an arbitrarily point in $\mathds{R}^N_{<0}$. Then
$\lim\limits_{\|u\|\rightarrow +\infty}f(u)=+\infty.$
\end{lem}
\pf Set $\Pi_{u_0}=\left\{u=(u_1,\cdots,u_N)^T\,|\,u_i\leq (u_0)_i,\,i=1,\cdots,N\right\}$, we claim
$$\lim_{\|u\|\rightarrow +\infty,\;u\in \Pi_{u_0}}f(u)=+\infty.$$
As $Hess f=L>0$, we know $f$ is strictly convex to the downwards in $\mathds{R}^N_{<0}$.
Because $\nabla f=K-K_0$ is an embedding, $\nabla f(u)=0$ has a unique solution $u=u_0$.
For every real half-line $l_u=\{u_t\triangleq u_0+t(u-u_0),\,t\geq 0\}$, $u\in\Pi_{u_0}$,
the function $t\mapsto \varphi(t)=f(u_t)$ is strictly convex to the downwards.
Furthermore, $\varphi(0)=\varphi'(0)=0$, $\varphi''(0)>0$,
which imply $\varphi(t)$ is strictly increasing and converges to $+\infty$ as $t\rightarrow \infty$.
Denote $h(t)\triangleq \inf\,\{f(x)\,|\,x\in \partial B(u_0,t)\cap \Pi_{u_0}\}$.
We want to prove $\lim\limits_{t\rightarrow+\infty}h(t)=+\infty$.
If not, since $h(t)$ is a strictly monotone increasing function, we may find $M>0$, such that $h(t)<M$ for all $t\geq0$.
For each $k\in\ \mathbb{N}$, select  $x_k \in \partial B(u_0,k)\cap \Pi_{u_0}$ such that $f(x_k)<M$.
Denote $\overset{\circ}{x_k}=u_0+\frac{x_k-u_0}{k}\in \partial B(u_0,1)\cap \Pi_{u_0}$.
$\{\overset{\circ}{x_k}\}$ must has a convergent subsequence, denoted as $\{\overset{\circ}{x_k}\}$ also.
Suppose $\overset{\circ}{x_k}\rightarrow x^*\in \partial B(u_0,1)\cap \Pi_{u_0}$.
Then there must exists a positive integer $a>0$, such that $f(x_t^*)>M$ for all $t\geq a$, where $x_t^*= u_0+t(x^*-u_0)$.
Select $\delta>0$ so that for all $x\in \mathscr{U}_a\triangleq B(x_a^*,\delta)\cap \partial B(u_0,a)$ we have $f(x)>M$. However, $\mathscr{U}_1\triangleq B(x^*,\,\delta/a)\cap \partial B(u_0,1)$ is an open neighborhood (relative to $\partial B(u_0,1)$) of $x^*$. So there exists at least one $m>a$ such that $\overset{\circ}{x_m}\in \mathscr{U}_1$. Then  $u_0+a(\overset{\circ}{x_m}-u_0)\in \mathscr{U}_a$, hence
$$f(x_m)=f(u_0+m(\overset{\circ}{x_m}-u_0))>f(u_0+a(\overset{\circ}{x_m}-u_0))>M,$$
which contradicts the selection of $x_k$, that is, $f(x_k)<M$ (for $\forall k$).
Thus we get the claim above.

For any sequence $\{u_n\}_{n=1}^{\infty}$ in $\mathds{R}^N_{<0}$ satisfying $\|u_n\|\rightarrow +\infty$, due to what we have proved above,
we may assume that all $u_n$ are in the set $\mathds{R}^N_{<0}-\Pi_{u_0}$.
Denote $\overset{\circ}{u_n}$ as the intersection point of half-line $\overrightarrow{uu_n}$ and $\partial B(u_0,1)$. Then $\overset{\circ}{u_n}$ has at least a accumulation point in $\partial B(u_0,1)\cap \partial\Pi_{u_0}$. Using similar explanation in the above paragraph, we can get $f(u_n)\rightarrow +\infty$. Thus we have $$\lim\limits_{\|u\|\rightarrow +\infty,\;u\in\mathds{R}^N_{<0}}f(u)=+\infty.$$
\qed

\begin{ac}
The first author would like to show his greatest respect to Professor Gang Tian. The first author would also like to thank Professor Guanxiang Wang, Professor Feng Luo, Dr. Xu Yiyan for helpful discussions. The second author would like to thank Professor Zhang Xiao for the invitation to AMSS. Both authors would like to give special thanks to Dr. Zhou Da, Wang Ding for their encouragements.
\end{ac}

\end{document}